\documentclass[12pt]{amsart}
\usepackage{mathrsfs, latexsym,amsmath,amsfonts,amssymb, cite,
graphicx, psfrag, verbatim}

\title{On the $p$-parts of quadratic Weyl group multiple Dirichlet
series}
\author{Gautam Chinta}
\author{Solomon Friedberg}
\author{Paul E. Gunnells}
 
\address{Department of Mathematics,
The City College of CUNY,
New York, NY 10031, USA
}
\email{chinta@sci.ccny.cuny.edu}
 
\address{Department of Mathematics,
Boston College,
Chestnut Hill, MA 02467, USA}
\email{friedber@bc.edu}

\address{Department of Mathematics and Statistics, University of
Massachusetts, Amherst, MA 01003, USA}
\email{gunnells@math.umass.edu}
 
\thanks{Research supported by NSF grants DMS--0354534 (Chinta), 
DMS-0353964 (Friedberg) (FRG grants), and DMS-0401525 (Gunnells).}

\date{\today}

\DeclareMathOperator{\sgn}{sgn}

\DeclareMathOperator{\Supp}{Supp}
\DeclareMathOperator{\adj}{adj.}

\newenvironment{proofof}[1]{\vspace*{.1in}\noindent{\em
Proof{#1}. \/}}{\qed\vspace{3ex}}
\begin{document}

\begin{abstract} 
Let $\Phi$ be a reduced root system of rank $r$.  A {\it Weyl group
multiple Dirichlet series} for $\Phi$ is a Dirichlet series in $r$
complex variables $s_1,\dots,s_r$, initially converging for $\Re(s_i)$
sufficiently large, which has meromorphic continuation to ${\mathbb
C}^r$ and satisfies functional equations under the transformations of
${\mathbb C}^r$ corresponding to the Weyl group of $\Phi$.  Two
constructions of such series are available, one
\cite{wmd1,wmd2,wmd3,wmd4} based on summing products of $n$-th order
Gauss sums, the second \cite{qmds} based on averaging a certain group
action over the Weyl group.  In each case, the essential work occurs
at a generic prime $p$; the local factors, satisfying local functional
equations, are then pieced into a global object.  In this paper we
study these constructions and the relationship between them.  First we
extend the averaging construction to obtain twisted Weyl group
multiple Dirichlet series, whose $p$-parts are given by evaluating
certain rational functions in $r$ variables.  Then we develop
properties of such a rational function, giving its precise
denominator, showing that the nonzero coefficients of its numerator
are indexed by points that are contained in a certain convex polytope,
determining the coefficients corresponding to the vertices, and
showing that in the untwisted case the rational function is uniquely
determined from its polar behavior and the local functional equations.
We also give evidence that in the case $\Phi=A_r$, the $p$-part
obtained here exactly matches the $p$-part of the twisted multiple
Dirichlet series introduced in \cite{wmd3} when $n=2$.
\end{abstract}


\maketitle

\newcommand{\Q}{\ensuremath{{\mathbb Q}}}
\newcommand{\R}{\ensuremath{{\mathbb R}}}
\newcommand{\A}{\ensuremath{{\mathbb A}}}
\newcommand{\C}{\ensuremath{{\mathbb C}}}
\newcommand{\Z}{\ensuremath{{\mathbb Z}}}
\newcommand{\x}{\ensuremath{{\bf { x}}}}
\renewcommand{\l}{\ensuremath{{\bf { l}}}}
\newcommand{\s}{\ensuremath{{\bf { s}}}}
\newcommand{\bbf}{\ensuremath{{\bf { f}}}}
\newcommand{\E}{{\mathcal E}}
\newcommand{\II}{\ensuremath{{\bf { I}}}}
\newcommand{\f}{\ensuremath{{\mathfrak { f}}}}
\newcommand{\fS}{\ensuremath{{\mathfrak { S}}}}
\newcommand{\sL}{\mathscr{L}}
\newcommand{\sR}{\mathscr{R}}
\newcommand{\sO}{\mathscr{O}}

\newcommand{\fixme}{\textbf{fixme.}}

\newcommand{\ttt}{\ensuremath{{\bf {t}}}}
\newcommand{\im}{\mbox{Im}}
\newcommand{\res}{\mathop{\mbox{Res}}}
\newcommand{\re}{\mbox{Re}}
\newcommand{\nin}{\noindent}
\newcommand{\mb}{\medbreak}
\newcommand{\adjacent}[2]{#1\sim#2}
\newcommand{\half}{{\smfrac{1}{2}}}
\newcommand{\G}{{\gamma}}
\newcommand{\B}{{\beta}}
\newcommand{\I}{{\mathcal I}}
\newcommand{\J}{{\mathcal J}}
\newcommand{\pr}{{\prime}}
\newcommand{\innprod}[2]{\langle #1, #2 \rangle}
\newcommand{\posroots}{\Phi^{+}}
\newcommand{\weight}{\varpi}
\newcommand{\norm}[1]{\left|\left|#1\right|\right|}
\newcommand{\abs}[1]{\left|#1\right|}
\newcommand{\dlim}[1]{\displaystyle {\lim_{#1}}}
\newcommand{\dsup}[1]{\displaystyle {\sup_{#1}}}
\newcommand{\dinf}[1]{\displaystyle {\inf_{#1}}}
\newcommand{\dmax}[1]{\displaystyle {\max_{#1}}}
\newcommand{\dmin}[1]{\displaystyle {\min_{#1}}}
\newcommand{\dprod}[2]{\displaystyle {\prod_{#1}^{#2}}}
\newcommand{\dcoprod}[2]{\displaystyle {\coprod_{#1}^{#2}}}
\newcommand{\dsum}[2]{\displaystyle {\sum_{#1}^{#2}\:}}
\newcommand{\prs}[2]{\left(\frac{#1}{#2}\right)}
\newcommand{\smfrac}[2]{{\textstyle \frac{#1}{#2}}} 
\newcommand{\nn}{\nonumber}
\newcommand{\no}{{\mathbb N}}
\newcommand{\cO}{{\mathcal O}}

\newtheorem{theorem}{Theorem}
\newtheorem{lemma}[theorem]{Lemma}
\newtheorem{prp}[theorem]{Proposition}
\newtheorem{proposition}[theorem]{Proposition}
\newtheorem{corollary}[theorem]{Corollary}
\newtheorem{conjecture}[theorem]{Conjecture}
\newtheorem{fact}[theorem]{Fact}
\newtheorem{claim}[theorem]{Claim}
\theoremstyle{definition}
\newtheorem{example}[theorem]{Example}
\newtheorem{remark}[theorem]{Remark}

\newtheorem{defin}[theorem]{Definition}

\numberwithin{theorem}{section}
\numberwithin{equation}{section}

\section{Introduction} Let $\Phi$ be a reduced root system of rank
$r$.  A {\it Weyl group multiple Dirichlet series} for $\Phi$ is a
Dirichlet series in $r$ complex variables $s_1,\cdots,s_r$, initially
converging for $\Re(s_i)$ sufficiently large, which has meromorphic
continuation to $\C^r$ and satisfies functional equations under the
transformations of $\C^r$ corresponding to the Weyl group of $\Phi$.
For example, every Langlands $L$-function $L(s,\pi,r)$ is a Dirichlet
series that is expected to have a functional equation of type $A_1$,
since $s\mapsto 1-s$ is a transformation of order $2$.  Let $n$ be a
fixed positive integer and let $K$ be a fixed global field containing
the $2n$-th roots of unity\footnote{That $K$ contain the $n$-th roots
of unity is essential.  The extra requirement that $K$ contain the
$2n$-th roots is made for convenience.}. In \cite{wmd1,wmd2,wmd3,wmd4}
a broad class of Weyl group multiple Dirichlet series was exhibited,
based on summing products of $n$-th order Gauss sums.  These series
are expected to be the Whittaker coefficients of minimal parabolic
Eisenstein series on the $n$-fold cover of the simply connected
algebraic group $G$ over $K$ whose root system is the dual of the root
system $\Phi$ (i.e.\ $\Phi$ is the root system of the
$L$-group~$^LG$).  In the articles cited above, the series was defined
for all $\Phi$ if $n$ is sufficiently large and for all $n$ if
$\Phi=A_r$, and its properties were established in the first of these
cases, and conjectured in the second.

Let $g(c)$ denotes the $n$-th order Gauss sum of modulus $c$ formed
from the $n$-th power residue symbol $\left(\frac{a}{c}\right)_n$.
Then the Chinese remainder theorem implies that for $(c_1,c_2)=1$ the
sum satisfies the twisted multiplicativity relation
$$g(c_1c_2)=\left(\frac{c_1}{c_2}\right)_n\left(\frac{c_2}{c_1}\right)_ng(c_1)\,g(c_2).$$
Similarly, the coefficients of the Weyl group multiple Dirichlet
series satisfy a twisted multiplicativity (see \cite{wmd3}, Eq.\ (2)).
Thus, though these series are {\it not} in general Euler products,
they may be reconstituted from a description of their $p$-part for a
given prime $p$ of norm $q.$ We denote the (arithmetic part of the)
coefficients at such a prime by $H(p^{k_1},\dots,p^{k_r})$.  (See
\cite{wmd2,wmd3,wmd4}.)

The Weyl group multiple Dirichlet series in the case $n=2$ are of
particular interest; for example, the study of automorphic forms on
double covers of classical groups is related to the theta
correspondence.  For $n=2$, an alternative method of defining a Weyl
group multiple Dirchlet series, valid for {\it all} simply-laced root
systems, was given in \cite{qmds}.  (Though the arguments of $n$-th
order Gauss sums vary for $n>2$, for $n=2$ they are essentially fixed,
and the construction does not make direct use of Gauss sums.)  In this
approach, the $p$-part of the series is constructed as a certain
average over the Weyl group under a natural but non-obvious action. By
construction, this series satisfies the requisite continuation and
local functional equations, and these lead to a global functional
equation.  However, the individual coefficients
$H(p^{k_1},\dots,p^{k_r})$ are not identified.  Nor is the size of the
space of $p$-parts satisfying the local functional equations apparent
(is the series, for example, uniquely determined by them?).  By
contrast, the approach of \cite{wmd3} provides a combinatorial
description for the coefficients of a series of type $A_r$ that is
expected to have continuation and functional equations.  However,
except in low-rank cases, these properties are only conjectured.
Moreover, it is not obvious that these two very different
constructions are even related.  In this paper we take a number of
steps towards addressing these issues.

To describe our results further, let us review some basic features of
these two constructions.  We begin with the series of
\cite{wmd1,wmd2,wmd3,wmd4}.  Let $\Phi$ be a simply-laced root system
with Weyl group $W$.  Fix a decomposition of $\Phi = \Phi^{+}\cup
\Phi^{-}$ into positive and negative roots, and let
$\alpha_1,\dots,\alpha_r$ be the simple roots.  Let
$\rho=\tfrac12\sum_{\alpha\in\Phi^+}\alpha$.  If $w\in W$ and
$\rho-w\rho=\sum_{i=1}^r k_i\alpha_i$, then the point
$(k_1,\dots,k_r)$ is called {\it stable}.  The stable points form the
vertices of a convex polytope called the \emph{permutahedron} attached
to $\Phi$.  The geometry of this object will be of great concern
below.  The series described in \cite{wmd1,wmd2,wmd3,wmd4} have the
properties that: (a) the coefficient $H(p^{k_1},\dots,p^{k_r})$
associated to the stable vertex attached to $w$ is a (specific)
product of $l(w)$ $n$-th order Gauss sums, where $l(w)$ is the length
of $w$; (b) if $n$ is sufficiently large (the {\it stable case}) then
all other coefficients are zero; (c) for all $n\geq1$, all
coefficients corresponding to points outside the permutahedron are
zero.

The series described in \cite{wmd1,wmd2,wmd3,wmd4} come in two
flavors: untwisted and twisted.  The untwisted series should
correspond to the first Whittaker coefficient for the Eisenstein
series; the twisted coefficients should correspond to higher Whittaker
coefficients.  The adjective ``twisted'' is used since when the
relevant indices are relatively prime, this corresponds to twisting
the original Weyl group multiple Dirichlet series by characters.  The
description of the previous paragraph is for the untwisted series;
however, it remains valid provided one does not twist by powers of
$p$.  It is precisely when one looks at the $p$-part of a series with
twisting by powers of $p$ that new features arise.  In this case, the
permutahedron expands (in a specific way; see \cite{wmd4}) to take
into account the twisting.  The vertices of this expanded
permutahedron are called the {\it twisted stable} points.  With the
permutahedron replaced by this expanded one and the stable vertices by
the twisted stable ones, statements (a), (b), and (c) above are still
expected to hold.  In \cite{wmd3} precise conjectures concerning the
interior coefficients are formulated for all $n$ in the case
$\Phi=A_r$, while in \cite{wmd4} the twisted series are defined for
all reduced root systems $\Phi$ and their continuations and functional
equations established, provided $n$ is sufficiently large (depending
on $\Phi$).  One expects the existence of twisted Weyl group multiple
Dirichlet series for all $\Phi$ and all $n$, and the formula of
\cite{wmd4} for the coefficients at the twisted stable vertices (once
again, a product of $l(w)$ Gauss sums; see
Theorem~\ref{thm:comparison} below) is expected to remain valid for
all $n$.  For $A_r$ this is consistent with the conjectures in
\cite{wmd3}; see \cite{wmd4}.

We turn to the series described in \cite{qmds}.  Suppose that the root
system $\Phi$ is simply-laced.  Let $\Lambda
= \Lambda_{\Phi}$ be the root lattice.  For $\lambda=\sum_{i=1}^r
k_i\alpha_i\in\Lambda$, let $d(\alpha)=\sum_{i=1}^r k_i$ be the usual
height function, and let $\x^{\lambda}$ be the monomial $x_{1}^{k_{1}}\dotsb
x_{r}^{k_{r}}$. Let $p$ be a prime of norm $q$.  In \cite{qmds}, the
authors construct a rational function $f (\x)$:
\[
f (\x) = \frac{N (\x)}{D (\x)}, \quad D (\x) := \prod_{\alpha \in
\posroots} (1-q^{d (\alpha)-1}\x^{2\alpha}),
\]
which is invariant under a certain $W$-action.  The construction will
be described in detail in the next section.  It is shown that
$f(q^{-s_1},\dots,q^{-s_r})$ is the $p$-part of a multiple Dirichlet
series (associated to $\Phi$ and $n=2$) which has meromorphic
continuation to ${\mathbb C}^r$ and a group of functional equations
isomorphic to $W$.
 
The work of \cite{qmds} concerns only untwisted series.  However, the
method can be adapted to construct twisted series as well.  We carry
out this extension, based on defining a {\it twisted action} of $W$ on
the field of rational functions in $x_1,\dots,x_r$, in Section 2
below.  Thus, our first result is a wider class of Weyl group multiple
Dirichlet series with continuation to $\C^r$.

With this extension in hand, we pursue two goals.  One is to develop
properties of the construction of \cite{qmds} and of this extension;
the second is to study its relation to the series of
\cite{wmd1,wmd2,wmd3,wmd4}.  The twisted rational function
$f(\x;\ell)$ (here $\ell\in(\Z_{\ge0})^r$ is an index that specifies
the twisting) is of the form $N(\x;\ell)/D(\x)$.  Write
\[
N (\x;\ell) = \sum_{\lambda \in \Lambda} a_{\lambda}\x^{\lambda}.
\]
Then we will show
\begin{itemize}
\item  $N(\x;\ell)$ is a polynomial (Theorem \ref{thm:NIsPoly}).
\item The support of $N(\x;\ell)$ is contained in a certain convex
polytope $\Pi_{\theta(\ell)}$ (Theorem \ref{thm:convex}).
\item The coefficients $a_\lambda$ for $\lambda$ a vertex of
$\Pi_{\theta(\ell)}$ coincide with the stable coefficients described
in \cite{wmd4} (Theorem \ref{thm:comparison}).
\item In the untwisted case, the series obtained in \cite{qmds} is the
{\it unique} series with denominator $D(\x)$ satisfying the desired
functional equations and having leading coefficient $1$ (Corollary
\ref{cor:untwistedunique}).
\end{itemize}
Even in the untwisted case each of these results is new.  In that case
$\Pi_{\theta(\ell)}$ is the convex hull in $\Lambda_{\R} = \Lambda
\otimes \R$ of the points
\[
\rho -w\rho, \quad w\in W.
\]

Based on these results, it seems natural to conjecture that if
$\Phi=A_r$ then for a given prime $p$ the series constructed here,
generalizing the construction of \cite{qmds}, is in fact {\it
identically the same} as the $p$-part of the series constructed using
Gelfand-Tsetlin patterns in \cite{wmd3}.  We note some computational
evidence for this as well (see the remarks after
Example~\ref{ex:A2polys}).

The remainder of this paper is organized as follows.  In Section 2, we
construct a twisted analogue of the quadratic Weyl group multiple
Dirichlet series of \cite{qmds}, valid for any simply laced root
system $\Phi$.  The $p$-part is obtained from a rational function
$f(\x;\ell)$.  We identify the denominator of this function, capturing
the polar behavior of the series.  We provide some explicit examples
when $\Phi=A_2$, note that in each case this series is the same as the
series of \cite{wmd3}, and describe some computational evidence that
suggests that this identification holds in general.  In Section 3, we
study the support of the numerator $N(\x;\ell)$, and show that it is
contained in a polytope $\Pi_{\theta(\ell)}$ as described above.  This
result depends on describing the polytope group-theoretically.  In
Section 4 we find the coefficients at the vertices of the polytope,
that is, the stable coefficients (both twisted and untwisted).  We
show that these coefficients are all determined from the degree $0$
coefficient of $N(\x;\ell)$.  As noted, these match the coefficients
of \cite{wmd3,wmd4}.  In Section 5, we discuss the remaining
coefficients of $N(\x;\ell)$.  In the untwisted case, we show that
these coefficients are all determined by the functional equations and
the stable coefficients, hence the series is uniquely determined by
the denominator, functional equations, and constant coefficient.
Surprisingly, in the twisted case, this uniqueness no longer holds,
and we explain why this is the case, due to the existence of other
regular orbits for the $W$-action on the points in
$\Pi_{\theta(\ell)}$.  Finally, in Section 6 we explain how to piece
together the $p$-parts described here into a twisted multiple
Dirichlet series with analytic continuation and group of functional
equations.

We thank Daniel Bump and Jim Humphreys for helpful conversations.  We
also thank Ben Brubaker and Daniel Bump for sharing software that
allowed us to compare our results with those of \cite{wmd3}.  Finally,
we thank the organizers of the Stanford Multiple Dirichlet Series
Workshop, held at Stanford University in the summer of 2006 and
supported by an NSF Focused Research Group grant, where substantial
portions of this work were carried out.


\section{Rationality and denominator}\label{s:ratandden}

We first develop the definition of the rational function $f(\x;\ell)$
which gives the $p$-part of a twisted quadratic Weyl group multiple
Dirichlet series.  We emphasize that this construction, based on
averaging over a suitable group action, was introduced and carried out
for the untwisted series in \cite{qmds}.  Here we modify the group
action in \cite{qmds} to accomodate twisting, and this allows us to
generalize those results using similar methods.  Once these $p$-parts
satisfy the desired local functional equations, it follows that they
may be combined by twisted multiplicativity to give a $W$-invariant
function on $\C^r$; see Section 6 below.

Let $F=\C({\bf x})=\C(x_1,x_2,\ldots, x_r)$ be the field of rational
functions in the variables $x_1,x_2,\ldots, x_r$.  
Let $ \ell=(l_1,
\ldots, l_r)$ denote an $r$-tuple of nonnegative integers.  This is
the {\em twisting parameter}. We define an $\ell$-twisted action of
$W$ on $F.$ This is done in stages.

Let $\sigma_{1},\dotsc ,\sigma_{r}\in W$ be the simple reflections
corresponding to the simple roots, and let $q$ be a fixed prime power.
First, for $\x=(x_1,x_2,\ldots, x_r)$ define $\sigma_i\x=\x^\prime$,
where
\begin{equation}\label{eqn:wiaction1}
x_j^\prime=\left\{\begin{array}{ll}
x_i x_j\sqrt{q} & \mbox{\ if $i$ and $j$ are adjacent,}\\
1/(q x_j) & \mbox{\ if $i=j$, and}\\
x_j & \mbox{\ otherwise.}
\end{array}\right.
\end{equation}

Next, define $\epsilon_i \x=\x^\prime$, where
\begin{equation}
x_j^\prime=\left\{\begin{array}{ll}
-x_j & \mbox{\ if $i$ and $j$ are adjacent,}\\
x_j & \mbox{\ otherwise.}
\end{array}\right.
\end{equation}
For $f\in F$ define
\begin{equation}
f_{i,\ell}^+(\x)=\frac{f(\x)+(-1)^{l_i}f(\epsilon_i \x)}{2}\mbox{ \ \ and\ \ }
 f_{i,\ell}^-(\x)=\frac{f(\x)-(-1)^{l_i}f(\epsilon_i \x)}{2}.
\end{equation}

Finally we can define the action of $W$ on $F$ for a generator
$\sigma_i\in W:$
\begin{equation}\label{wiactionTwisted}
(f|_\ell\sigma_i)(\x)=-\frac{1-qx_i}{qx_i(1-x_i)}(x_i\sqrt{q})^{l_i}
f_{i,\ell}^+(\sigma_i\x)+ 
\frac{1}{x_i\sqrt{q}}(x_i\sqrt{q})^{l_i}f_{i,\ell}^-(\sigma_i\x).
\end{equation}

Note that despite the presence of the $\sqrt{q}$'s, $(f|_\ell\sigma_i)(\x)$
will in fact be a power series in $q.$  For example,  suppose $l_i$ is
even.  Then in the function $f,$ each of the neighbors $x_j$ of $x_i$ gets 
replaced by $x_ix_j\sqrt{q}.$  But the sum of  degrees of the neighbors
of $x_i$ is even in $f_{i,\ell}^+$ and odd in $f_{i,\ell}^-.$
Therefore both summands in \eqref{wiactionTwisted} are in $K.$  
The argument is similar for $l_i$ odd.

For  $\ell = (0,\dotsc ,0)$, the action
\eqref{wiactionTwisted} coincides with the action of $\sigma_{i}$ on
$F$ from \cite[(3.13)]{qmds}.  As in \cite{qmds}, one can prove the
following lemma; the proof is essentially the same as that of
\cite[Lemma 3.2]{qmds}, and we omit the details.

\begin{lemma}
The definition \eqref{wiactionTwisted} extends to give an action of
$W$ on $F$.
\end{lemma}

The main result of Section 3 of \cite{qmds} is the construction of a
$W$-invariant rational function with certain limiting behavior.  This
is summarized below.

\begin{theorem}\label{thm:tw} Define
$$\Delta(\x)=\prod_{\alpha\in \Phi^+}(1-q^{d(\alpha)}\x^{2\alpha}),$$
$$j(w,\x)=\Delta(\x)/\Delta(w\x)=\sgn(w)q^{d(\rho -w^{-1}\rho )}\x^{2 (\rho -w^{-1}\rho )},$$
and
$$f(\x;\ell)=\Delta(\x)^{-1}\sum_{w\in W} j(w,\x) (1|_\ell w)(\x).$$
Then $f(\x;\ell)$ is a $W$-invariant rational function satisfying
\begin{enumerate}
\item for each $i=1,2,\ldots, r$, the function $f$  satisfies the following
  limiting condition: if $x_j= 0$ for every $j$ adjacent to $i$, then
\begin{equation}
f(\x;\ell)(1-x_i)^{m_i} \mbox{\ is independent of \ } x_i,
\end{equation}
where $m_i$ is 0 if $l_i$ is even and is 1 otherwise.
\item $f(0,0,\ldots, 0;\ell)=1.$
\end{enumerate}
\end{theorem}

In \cite{qmds} this theorem is proven only for the untwisted action,
i.e.~only for $\ell=(0,\ldots,0)$, but essentially the same proof
will work here.  We will not repeat the argument.  We will also need
below the analogues for the twisted action of Lemma 3.3(c) and Lemma
3.9 of \cite{qmds}.  We state these without proof.

\begin{lemma}\label{lemma:analogues}
Let $w\in W$ and $\ell$ be an $r$-tuple of nonnegative integers.
\begin{enumerate}
\item[(a)] Let $g,h\in F.$ 
If $g$ is an even function of all the $x_j$, then 
$$(gh|_\ell w)(\x)=g(w\x)\cdot (h|_\ell w)(\x).$$
\item[(b)] $\x^{\rho-w\rho}(1|_\ell w)(\x)$ is regular at the origin.
\item[(c)] $\x^{\rho-\sigma_iw\rho}\left(\smfrac{1}{x_i}|_\ell w\right)(\x)$ is
  regular at the origin for $i=1,2,\ldots, r$.
\end{enumerate}
\end{lemma}

Put \[
f (\x;\ell) = \frac{N (\x;\ell)}{D (\x)}, \quad D (\x) := \prod_{\alpha \in
\Phi^{+}} (1-q^{d (\alpha)-1}\x^{2\alpha}).
\]

\begin{theorem}\label{thm:NIsPoly}
$N(\x;\ell)$ is a polynomial.
\end{theorem}

\begin{proof}
First, we will show that $f(\x;\ell)D(\x)\Delta(\x)$ is a polynomial.  
In fact, we will prove, for any $w\in W,$
\begin{equation}\label{eqn:Dw}
j(w,\x)(1|_\ell w)(\x)D_w(\x)
\end{equation}
is a polynomial, where
$$D_w(\x)=\prod_{\substack{\alpha\in\Phi(w)}} 
(1-q^{d (\alpha)-1}\x^{2\alpha}).
$$
Here $\Phi(w)$ denotes the subset of positive roots made negative
by $w$.  The proof that \eqref{eqn:Dw} is a polynomial will be by
induction on the length $l(w)$ of $w.$ When $w$ is the identity, there
is nothing to prove.  Suppose that for $w\in W$, \eqref{eqn:Dw} is a
polynomial,
$$j(w,\x)(1|_\ell w)(\x)D_w(\x)=P(\x),  $$
say.
 Let $\sigma_i$ be a simple
reflection such that $l(w\sigma_i)=l(w) +1.$  Then 
\begin{eqnarray}\nonumber
\lefteqn{j(w\sigma_i,\x)(1|_\ell w\sigma_i)(\x)D_{w\sigma_i}(\x)}\\
\nonumber
&=& j(\sigma_i,\x)\left(\left.
\frac {P}{D_w}\right|_\ell \sigma_i\right)(\x)\cdot
D_{w\sigma_i}(\x)  \\ 
\label{eqn:jpd}
&=& \frac{ j(\sigma_i,\x)(P|_\ell \sigma_i)(\x)}
{D_w(\sigma_i\x)}\cdot
D_{w\sigma_i}(\x), \text{\ \ \ \ by Lemma 2.3(a)}.
\end{eqnarray}
By the definition \eqref{wiactionTwisted} of the action of $\sigma_i,$
we can write $(P|_\ell\sigma_i)(\x)$ as $P_2(\x)/(1-x_i^2)$ where
$P_2$ is a Laurent polynomial in the $x_i.$ However, as Lemma
\ref{lemma:analogues} (b) implies that \eqref{eqn:Dw} is regular at
the origin, it follows that $P_2$ is a polynomial.  Moreover, the
denominator $D_w(\sigma_i\x)$ is equal to
$$\prod_{\substack{\alpha\in\Phi(w)}} 
(1-q^{d(\sigma_i\alpha)-1}\x^{2\sigma_i\alpha})=
D_{w\sigma_i}(\x)/(1-x_i^2)$$
Plugging this back into \eqref{eqn:jpd}, we conclude that
$$j(w\sigma_i,\x)(1|_\ell w\sigma_i)(\x)D_{w\sigma_i}(\x)=
j(\sigma_i, \x)P_2(\x)
$$
is polynomial.  Therefore $f(\x)D(\x)\Delta(\x)=N(\x)\Delta(\x)$ is a
polynomial. 

To complete the proof, it will suffice to show that
\begin{equation}
 h(\x)= \sum_{w\in W} j(w,\x) (1|_\ell w)(\x)
\end{equation} is
divisible by $\Delta(\x).$   This function satisfies 
$$h( \x)=j(\sigma_i, \x) (h|_\ell \sigma_i)(\x).$$
Setting $x_i=\pm 1/\sqrt{q}$ in the above equation, we deduce that
$$h(\x)\bigr|_{x_i=\pm 1/\sqrt{q}}=- h(\x)\bigr|_{x_i=\pm 1/\sqrt{q}},$$
and that hence $1-qx_i^2$ divides $h(\x).$  Similarly, we can show that
$1-q^{d(\alpha)}\x^{2\alpha} $ divides $h(\x)$ for all positive roots
$\alpha.$ This completes the proof of the theorem.
\end{proof}

\begin{remark}
Though we do not require this fact, we point out that that the proofs
of Lemma \ref{lemma:analogues} and Theorem \ref{thm:NIsPoly} actually
show that $N$ is a polynomial in $q$ as well.
\end{remark}

\begin{example}\label{ex:A2polys}
Here we give some examples of the polynomials $N (\x ; \ell)$ for the
root system $A_{2}$.  We use $(x,y)$ for the variables instead of
$(x_{1},x_{2})$.  To get the rational function $f (\x ; \ell)$, each
of the following should be divided by $(1-x^{2}) (1-y^{2}) (1-qx^{2}y^{2})$.

\begin{itemize}
\item $\ell = (0,0)$: $N = 1 + x + y - x^2 y - x y^2 - x^2 y^2$.\vskip 3pt
\item $\ell = (1,0)$: $N =1-x^2+y+(q-1) x^2 y+q x^3 y-q x^2 y^3-q x^3 y^3$.\vskip 3pt
\item $\ell = (1,1)$: $N =1-x^2-y^2+(1-q) x^2 y^2+q x^4 y^2+q x^2 y^4-q x^4 y^4$.\vskip 3pt
\item $\ell = (2,0)$: $N =1+(q-1) x^2+q x^3+y+(q-1) x^2 y-q x^4 y+\left(q^2-q\right) x^2
   y^2+\left(q^{2}-q\right) x^3 y^2+\left(q^{2}-q\right) x^2
   y^3+\left(q-q^2\right) x^4 y^3-q^2 x^3 y^4-q^2 x^4 y^4 $.\vskip 3pt
\item $\ell = (2,1)$: $N =1+(q-1) x^2+q x^3-y^2+\left(1-2 q+q^2\right) x^2
   y^2+\left(q^2-q\right) x^3 y^2+\left(q-q^2\right) x^4 y^2-q^2
   x^5 y^2+\left(q-q^2\right) x^2 y^4+\left(q^{3}-q^2\right) x^3
   y^4+\left(q^{2}-q\right) x^4 y^4+$\\
   $\left(q^2-q^3\right) x^5 y^4+q^3
   x^3 y^5-q^3 x^5 y^5 $.\vskip 3pt
\item $\ell = (2,2)$: $1+(q-1) x^2+q x^3+(q-1) y^2+\left(1-3 q+2 q^2\right) x^2
   y^2+\left(q^2-q\right) x^3 y^2+\left(q-2 q^2+q^3\right) x^4
   y^2+\left(q^{3}-q^2\right) x^5 y^2+q y^3+\left(q^{2}-q\right) x^2
   y^3+\left(q^{3}-q^2\right) x^4 y^3-q^3 x^6 y^3+\left(q-2
   q^2+q^3\right) x^2 y^4+\left(q^{3}-q^2\right) x^3 y^4+\left(q^4-2 q^3+2
   q^2-q\right) x^4 y^4+$\\
   $\left(q^2-2 q^3+q^4\right) x^5
   y^4+\left(q^{3}-q^2\right) x^2 y^5+\left(q^2-2 q^3+q^4\right) x^4
   y^5+\left(q^3-q^4\right) x^6 y^5-q^3 x^3 y^6+\left(q^3-q^4\right)
   x^5 y^6-q^4 x^6 y^6$.
\end{itemize}

\end{example}
The invariance of $f(\x;\ell)$ and the limiting conditions imply that the
Taylor series coefficients of $f$ can be used to construct the
$p$-part of a twisted multiple Dirichlet series with continuation to $\C^r$
and functional equation.   The $p$-parts are combined using twisted
multiplicativity, as explained in Section 6.

In fact, {\it we conjecture that for $\Phi = A_{r}$ this $p$-part will
coincide with the $p$-part coming from the construction of \cite{wmd3}
via Gelfand--Tsetlin patterns.\ } As evidence, we have computed
$N(\x;\ell)$ for the root systems $A_r$, $r\leq 4$, for all twisting
parameters $\ell$ with $\sum l_i\leq 6$.  For each of these 2597
series the polynomial $N (\x ; \ell)$ matches the numerator polynomial
predicted by the Gelfand--Tsetlin conjecture.


\section{The support of the numerator}

Recall that  $\Lambda$ denotes the root lattice of $\Phi$, and that for
$\lambda =\sum_{i=1}^r k_i\alpha_i\in
\Lambda$,  $\x^{\lambda}$ denotes the monomial $\prod_{i=1}^r x_i^{k_i}$.  Write 
\begin{equation}\label{eq:expressionforN}
N (\x ; \ell ) = \sum_{\lambda \in \Lambda}
a_{\lambda}\x^{\lambda },
\end{equation}
where $N (\x ;\ell)$ is the polynomial from Theorem \ref{thm:NIsPoly}.
Our goal in this section is to investigate the
support of $N (\x ;\ell)$, that is the set 
\[
\Supp N (\x ;\ell) = \{ \lambda \in \Lambda\mid a_{\lambda}\not =0 \}.
\] 
The main result (Theorem \ref{thm:convex}) is that $\Supp N (\x ;
\ell)$ is contained in a certain translated weight polytope for
$\Phi$.

We begin by discussing relations that the coefficients $a_\lambda$
must satisfy.
Recall that $N (\x ; \ell)$ is the numerator of the rational function
\[
f (\x) = \frac{N (\x; \ell)}{D (\x)}, \quad D (\x) := \prod_{\alpha \in
\posroots} (1-q^{d (\alpha)-1}\x^{2\alpha}),
\]
and that $f (\x)$ is invariant under a certain action of the Weyl
group $W$.
The $W$-invariance of $f (\x )$ and the simple form of the denominator $D (\x)$
imply that the numerator $N (\x; \ell )$
satisfies certain relations under the $W$-action.  In particular,
applying the reflection $\sigma_{k}$ gives the relation 
\begin{multline}\label{eq:actiononN}
(q^{2}\x^{2\alpha_{k}}-1)N (\x; \ell ) = q^{(1+l_{k}/2)}\x^{(l_{k}+1)\alpha_{k}} (1+\x^{\alpha_{k}})
(q\x^{\alpha_{k}}-1)N_{k,\ell}^{+} (\sigma_{k}\x)\\
+ q^{(3/2+l_{k}/2)}\x^{(l_{k}+1)\alpha_{k}}
(1-\x^{2\alpha_{k}})N_{k,\ell}^{-} (\sigma_{k} \x).
\end{multline}
Inserting \eqref{eq:expressionforN} in \eqref{eq:actiononN} and collecting terms, we find
\begin{equation}\label{eq:coeffs}
q^{2}a_{\lambda -2\alpha_{k}} - a_{\lambda} = \begin{cases}
q^{d (\beta)/2} (-a_{\mu }+ (1-1/q)a_{\mu+\alpha_{k}} +a_{\mu +2\alpha_{k}}/q)&\text{(even)},\\
\sqrt{q}\cdot q^{d (\beta)/2} (a_{\mu}-a_{\mu +2\alpha_{k}}/q^{2})&\text{(odd)},
\end{cases}
\end{equation}
In \eqref{eq:coeffs} 
we have used the notation
\begin{itemize}
\item $\mu = \sigma_{k}\lambda + (l_{k}+1)\alpha_{k}$,
\item $\beta =
\lambda -\sigma_{k}\lambda -l_{k}\alpha_{k} = \lambda -\mu +\alpha_{k}$, and 
\item $d \colon
\Lambda \rightarrow \Z$ is the usual height function on the root lattice;
\end{itemize}
we also are assuming that $\lambda \leq \mu$, in other words that the
difference $\mu -\lambda$ is a nonnegative sum of simple roots.
``Even/odd'' in \eqref{eq:coeffs} refers to the following.  Define $\varphi_{k}\colon
\Lambda\rightarrow \Z$ by 
\begin{equation}\label{eq:phikdef}
\varphi_{k} (\mu) = l_{k}+\sum_{\adjacent{j}{k}}k_{j},
\end{equation}
where $\mu =\sum k_{j}\alpha_{j}$, and where $\adjacent{j}{k}$ means
that the nodes labelled $j$ and $k$ are adjacent in the Dynkin diagram
for $\Phi$.
Then we are in the even/odd case according to whether $\varphi_{k}
(\mu)$ is even or odd.  Note that $\varphi_{k} (\mu) = \varphi_{k}
(\lambda)$, since $\mu -\lambda$ is a multiple of $\alpha_{k}$.

The functional equation (\ref{eq:coeffs}) can be simplified in the
even case.  Keeping the notation $\lambda,\mu,\beta$ as above, let us
apply the functional equation (\ref{eq:coeffs}) with $\lambda$
replaced by $\mu+2\alpha_k$.  Since
$\sigma_k(\mu+2\alpha_k)+\alpha_k=\lambda-2\alpha_k$ and
$\mu+2\alpha_k-(\lambda-2\alpha_k)+\alpha_k=-\beta+6\alpha_k,$ we
obtain
\begin{equation}\label{eq:coeffs-1}
q^2 a_{\mu} - a_{\mu+2\alpha_k} = \begin{cases}
q^{3-d (\beta)/2} (-a_{\lambda-2\alpha_k}+
(1-1/q)a_{\lambda-\alpha_{k}} +a_{\lambda}/q)&\text{(even)},\\
\sqrt{q}\cdot q^{3-d (\beta)/2} (a_{\lambda-2\alpha_k}-a_{\lambda}/q^{2})&\text{(odd)}.
\end{cases}
\end{equation}
Notice that $\mu$, $\lambda+2\alpha_k$ are identically equal except in the $k$-th position,
so $\varphi_k(\mu)$ is even if and only if $\varphi_k(\lambda+2\alpha_k)$ is even
(and if and only if $\varphi_k(\lambda)$ is even).
In the odd case, equations (\ref{eq:coeffs}), (\ref{eq:coeffs-1})
are equivalent.  However, in the even case, they are not; instead, it is
is easy to combine them to obtain
\begin{equation}\label{eq:neweven}
a_{\lambda}+qa_{\lambda-\alpha_k} =
q^{d (\beta)/2-1} (qa_{\mu }+a_{\mu+\alpha_{k}})
\qquad{\text{(even).}}
\end{equation}
Hence (\ref{eq:neweven}) holds for all $\lambda$ in the even case.
Conversely, in the even case (\ref{eq:neweven}) for all $\lambda$ implies (\ref{eq:coeffs}):
take $q$ times (\ref{eq:neweven}) with $\lambda$
replaced by $\lambda-\alpha_k$ and subtract the original  (\ref{eq:neweven}).
We summarize the above discussion in the following proposition:

\begin{proposition}\label{prop:coeff.fes}
Let $\lambda ,\mu \in \Lambda$ satisfy $\mu = \sigma_{k}\lambda +
(l_{k}+1)\alpha_{k}$, and suppose $\mu \geq \lambda$.  Put $\beta
=\lambda -\mu +\alpha_{k}$, and define
$\varphi_{k}\colon \Lambda \rightarrow \Z$ as in \eqref{eq:phikdef}.  Then the
coefficients of $N (\x ; \ell)$ satisfy 
\begin{subequations}
\begin{align}
qa_{\lambda-\alpha_k} + a_{\lambda} &= q^{d (\beta)/2} (a_{\mu
}+a_{\mu+\alpha_{k}}/q), \quad \text{if $\varphi_{k} (\mu)$ is even, and}\label{eq:coeffs.a}\\
q^{2}a_{\lambda -2\alpha_{k}} - a_{\lambda} &= \sqrt{q}\cdot q^{d
(\beta)/2} (a_{\mu}-a_{\mu +2\alpha_{k}}/q^{2}), \quad \text{if
$\varphi_{k} (\mu)$ is odd.}\label{eq:coeffs.b}
\end{align}
\end{subequations}
\end{proposition}

Now let $\varpi_{1},\dotsc ,\varpi_{r}$ be the fundamental weights of
the root system $\Phi$; since $\Phi$ is simply-laced, these form the
dual basis to the simple roots with respect to the $W$-invariant
scalar product $\innprod{\cdot}{\cdot}$.   Recall that the closure of the \emph{dominant chamber}
in $\Lambda_{\R} = \Lambda \otimes \R$ is defined by the inequalities 
\begin{equation}\label{eq:antidomineq}
\innprod{\varpi_{k}}{\x}\geq  0, \quad k=1,\dotsc ,r;
\end{equation}
here we have slightly abused notation to let
$\x$ denote a point in $\Lambda_{\R}$.  Given a twisting parameter
$\ell = (l_{1},\dotsc ,l_{r})$, let $\theta = \theta (\ell)$ be the
dominant weight 
\[
\theta  = \rho + \sum l_{k}\varpi_{k},
\]
where $\rho$ is the sum of the fundamental weights.  Note that
$\theta$ is \emph{regular}, that is the inequalities
\eqref{eq:antidomineq} are strict when evaluated on $\theta$.

Let $\Pi = \Pi_{\theta}$ be the convex hull in $\Lambda_{\R}$ of
the points 
\[
\theta - w\theta, \quad w\in W.
\]
Our goal is to prove the following theorem:
\begin{theorem}\label{thm:convex}
The support $\Supp N (\x ; \ell)$ is contained in the polytope $\Pi$.
\end{theorem}

\begin{example}
Figure \ref{fig:22support} shows the support for the polynomial $N
(x,y; (2,2))$ from Example \ref{ex:A2polys}.  The shaded hexagon is the
polygon $\Pi$.  The grey dots are nonzero coefficients of $N
(x,y; (2,2))$.

\begin{figure}[htb]
\psfrag{x}{$x$}
\psfrag{y}{$y$}
\begin{center}
\includegraphics[scale=0.5]{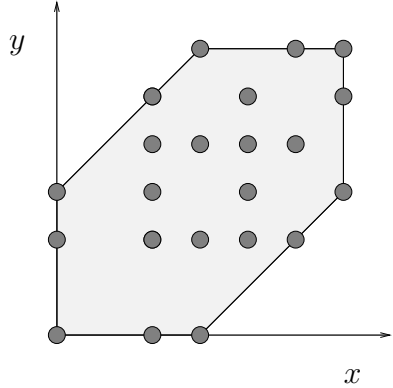}
\end{center}
\caption{The support of $N (x,y; (2,2))$\label{fig:22support}.}
\end{figure}
\end{example}

Before we prove the theorem we make a few remarks about the polytope
$\Pi$ and the geometry of the action
\eqref{eq:coeffs.a}, \eqref{eq:coeffs.b}.  After shifting $\Pi$ by
$\theta$, we see that $\Pi$ is isomorphic to the convex hull $\Pi '$
of the points $\{-w\theta \mid w\in W\}$, and is thus a
\emph{Coxeterhedron} or \emph{permutahedron} of type $W$.  In
particular, from the theory of such polytopes we know that the
vertices of $\Pi$ are exactly the points $\theta -w\theta$, and hence
that $\Pi$ has $|W|$ vertices.  Shifting makes the connection between
$\lambda$ and $\mu$ more apparent.  Suppose that $\lambda '$
(respectively, $\mu '$) is the vertex of $\Pi '$ obtained by
translating $\lambda$ (resp., $\mu$). Then $\lambda '$ and $\mu '$ are
related by $\sigma_{k}\lambda ' = \mu '$; in other words, after
shifting $\Pi$ to $\Pi '$, the functional equations
\eqref{eq:coeffs.a}, \eqref{eq:coeffs.b} relate coefficients of
monomials attached to weights that are connected by the usual
reflection action of $W$.

It is not hard describe a set of inequalities defining $\Pi$: it is
cut out by the system
\begin{equation}\label{eq:inequalities}
\innprod{w\weight_{i}}{\x - (\theta -w\theta)}\geq 0, \quad w\in W,
\quad i=1,\dotsc ,r.
\end{equation}
To prove this, one observes that the inequalities 
\[
\innprod{\weight_{i}}{\x}\geq 0, \quad i=1,\dotsc ,r
\]
define the facets containing the origin, and uses the fact that the
Weyl group $W$ acts transitively on the vertices by affine
transformations.  The same computation shows that the inequalities
labelled by $w$ are active at the vertex $\theta - w\theta$.

The system \eqref{eq:inequalities} is redundant.  We clarify this in
the following lemma, whose statement requires the notion of the
\emph{right descent set} of an element $w\in W$.  By definition, this
is the set $\sR (w) = \{\sigma_{i}\mid l (w\sigma_{i})<l (w) \}$.

\begin{lemma}\label{lem:redundant}
Let $\sigma_{j}\in \sR (w)$ and let $u = w\sigma_{j}$.  Then if $j\not
=k$, the inequalities
\[
\innprod{w\weight_{k}}{\x - (\theta - w\theta)}\geq 0
\]
and 
\[
\innprod{u\weight_{k}}{\x - (\theta - u\theta) }\geq 0
\]
are equivalent.
\end{lemma}

\begin{proof}
This follows since $\sigma_{j}\weight_{k} = \weight_{k}$ if $j\not
=k$.  Indeed, starting with the second inequality, we have 
\begin{align*}
\innprod{u\weight_{k}}{\x - (\theta - u\theta )} &=
\innprod{u\sigma_{j}\weight_{k}}{\x - (\theta -u\theta )} \\
&=
\innprod{u\sigma_{j}\weight_{k}}{\x - (\theta -u\theta ) +
w\theta  -w\theta  } \\
&= \innprod{w\weight_{k}}{\x - (\theta - w\theta)} +
\innprod{u\weight_{k}}{u\theta - w\theta)}.
\end{align*}
(In the last line we again used $\sigma_{j}\weight_{k} = \weight_{k}$
if $j\not =k$.)  By $W$-invariance, the second term on the last line
is
\[
\innprod{\weight_{k}}{\theta -\sigma_{j}\theta },
\]
which vanishes since $\theta -\sigma_{j}\theta$ is a multiple of $\alpha_{j}$.
This completes the proof.
\end{proof}

We will also need the following geometric lemma about $\Pi$, whose
statement uses the 
\emph{left descent set} $\sL (w)$ of an element $w\in W$.  By
definition $\sL (w) = \{\sigma_{i}\mid l (\sigma_{i}w)<l (w) \}$.
Recall that $\Phi (w)$ denotes the subset of the positive roots made
negative by $w$.  If $w=u\sigma_{k}$ and $l (w)=l (u)+1$, then from
the theory of Coxeter groups \cite{humph}
\begin{subequations}
\begin{align}
\Phi (w) &=\sigma_k\Phi(u)\cup \{\alpha_k \},\label{eqn:negrootscomputation}\\
\Phi(w^{-1})&=\Phi(u^{-1})\cup \{u\alpha_k \}.\label{eqn:negroots2}
\end{align}
\end{subequations}
\begin{lemma}\label{lem:string}
Let $\mu  = \theta - w\theta$ be a vertex of $\Pi$, and suppose $\sigma_{k}\in \sL
(w)$.  Then any lattice point of the form $\mu+m\alpha_{k}$,
where $m$ is a positive integer, lies outside $\Pi$.  Similarly, let $u
= \sigma_{k}w$ and let $\lambda = \theta - u\theta$.  Then any point of
the form $\lambda - m\alpha_{k}$, where $m$ is a positive integer,
lies outside $\Pi$.
\end{lemma}

\begin{proof}
For the first statement, 
it suffices to show that $\mu +m\alpha_{k}$ violates the inequalities
active at $\mu$, which are given by 
\[
\innprod{w\weight_{i}}{\x - (\theta -w\theta)}\geq 0, \quad i=1,\dotsc ,r.
\]
Thus we want to show 
\[
\innprod{w\weight_{i}}{\mu +m\alpha_{k} - \mu } =
\innprod{w\weight_{i}}{m\alpha_{k}} =
-m\innprod{u\weight_{i}}{\alpha_{k}} <0
\]
for at least one $i$.  In fact, we will prove the stronger statement
that 
\[
-m\innprod{u\weight_{i}}{\alpha_{k}}<0\quad \text{for all $i=1,\dotsc ,r$.}
\]
Since $m>0$, we must show 
\[
\innprod{u\weight_{i}}{\alpha_{k}} =
\innprod{\weight_{i}}{u^{-1}\alpha_{k}}>0, \quad i=1,\dotsc ,r. 
\]
This follows if and only if $\alpha_{k}\not \in \Phi (u^{-1})$.

So suppose on the contrary that $\alpha_{k}\in \Phi (u^{-1})$.  By
\eqref{eqn:negrootscomputation}, we
have 
\[
\Phi (w^{-1}) = \Phi(u^{-1}\sigma_k)=\sigma_k\Phi(u^{-1})\cup
\{\alpha_k \}.
\]
In particular $\sigma_{k} \Phi (u^{-1})$ consists of positive roots.
But if $\alpha_{k}\in \Phi (u^{-1})$ then $\sigma_{k}\alpha_{k}$ is
negative.  This contradiction completes the proof of the first
statement; the second statement is proved in almost exactly the same
way.
\end{proof}

\begin{figure}[htb]
\psfrag{l}{$\lambda $}
\psfrag{m}{$\mu $}
\psfrag{i}{$\text{id}$}
\psfrag{1}{$1$}
\psfrag{2}{$2$}
\psfrag{12}{$12$}
\psfrag{21}{$21$}
\psfrag{121}{$121$}
\psfrag{s1}{$\sigma_{1}$}
\begin{center}
\includegraphics[scale=0.4]{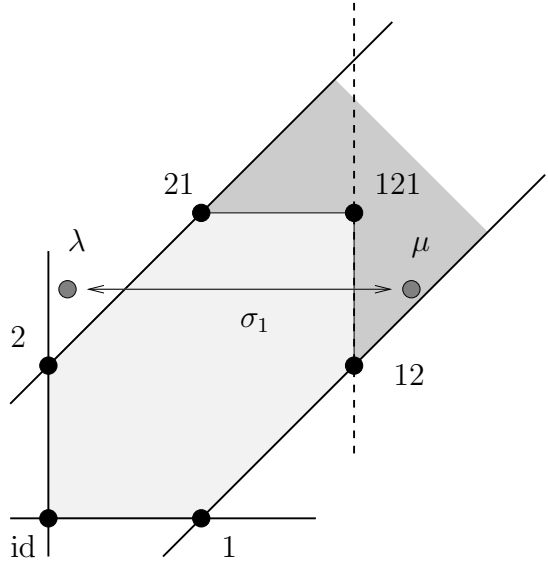}
\end{center}
\caption{Checking the inequalities at the vertex labelled $12$.  By
induction we assume that all desired inequalities hold at vertices
labelled by $w$ with $l (w)\leq 1$.  The point $\mu$ violating the
dashed inequality leads to a point $\lambda$ violating the
inequalities at $2 = \sigma_{1}\cdot 12$.\label{fig:inequalities}}
\end{figure}

\begin{proofof}{ of Theorem \ref{thm:convex}}
We use induction on the length.  Since $N (\x; \ell )$ is a polynomial, we
know that $a_{\lambda}=0$ if $\lambda$ violates the inequalities
active at the origin.  Indeed, otherwise $N (\x; \ell )$ would have polar terms.

Now let $w\in W$ satisfy $l (w)>0$, and suppose we have verified the
inequalities at all vertices $\theta -u\theta$ where $l (u)<l (w)$.  Let
$\sigma_{k}\in \sR (w)$, which is nonempty since $l
(w)>0$.  Then Lemma \ref{lem:redundant} implies that $a_{\lambda}=0$
unless $\lambda$ satisfies the inequalities 
\[
\innprod{w\weight_{j}}{\x - (\theta - w\theta)}\geq 0,
\]
for any $j\not =k$.  If $|\sR (w)|>1$, this shows that in fact all
desired inequalities hold for the support of $N (\x; \ell )$ at the
vertex $\theta - w\theta$.

Thus we assume $\sR (w) = \{\sigma_{k} \}$.  We must show
$a_{\lambda}=0$ if $\lambda$ violates 
\begin{equation}\label{eq:violate}
\innprod{w\weight_{k}}{\x - (\theta - w\theta)}\geq 0.
\end{equation}
Let $\sigma_{j}\in \sL (w)$; again $\sL (w)\not =\emptyset$ if $l
(w)>0$.  Choose $\mu \in \Lambda$ such that
\begin{enumerate}
\item $\mu$ violates \eqref{eq:violate}, 
\item $a_{\mu}\not =0$, and 
\item $a_{\mu'} = 0$ for all $\mu  '=\mu  +m\alpha_{j}$ with $m>0$.
\end{enumerate}
By Lemma \ref{lem:string} it is possible to find such a $\mu $.
Indeed the proof of Lemma \ref{lem:string} shows that if $\mu $
violates \eqref{eq:violate}, so do all the points $\mu 
+m\alpha_{j}$, $m>0$.  Since $N (\x;\ell )$ has bounded support there
must be final point in the support of $N (\x;\ell )$ on the ray $\mu
+m\alpha_{j}$. 

Now apply the relation \eqref{eq:coeffs} with $\sigma_{j}$ to $a_{\mu
}$, where $a_{\mu}$ is the first coefficient on the right side of
\eqref{eq:coeffs.a}, \eqref{eq:coeffs.b}. Note that since $a_{\mu +
m\alpha_{j}} = 0$ for $m>0$, as far as the right hand sides of these
equations are concerned it doesn't matter whether we are in the even
or the odd case.  Applying $\sigma_{j}$ produces the left hand side
$qa_{\lambda -\alpha_{j}}+a_{\lambda}$ in the even case, and $q^{2}a_{\lambda
-2\alpha_{j}}-a_\lambda$ and the odd case.  It is easy to see that
$a_{\lambda}$ and hence (Lemma \ref{lem:string}) $a_{\lambda
-\alpha_{j}}$ and $a_{\lambda -2\alpha_{j}}$ vanish by the induction
hypothesis, since $\lambda$ violates the inequalities active at $\theta 
-\sigma_{j}w\theta  $ (cf.~Figure \ref{fig:inequalities}).  Hence $a_{\mu}$ vanishes.  This shows that
$a_{\mu}=0$ unless $\mu$ satisfies \eqref{eq:violate}, and completes
the proof of the theorem.
\end{proofof}

\begin{remark}\label{rem:notnegative}
Theorem \ref{thm:convex} shows that the support of $N (\x ; \ell)$ is
contained in a translated weight polytope $\Pi$.  We caution the
reader that although $N (\x ; \ell)$ is constructed using the regular dominant
weight $\theta$, the polytope $\Pi $ is \emph{not} a translate of the
weight polytope $P$ for the representation with highest weight
$\theta$.  In fact, $\Pi$ is a translate of the weight polytope
attached to the representation with \emph{lowest weight} $-\theta$.
This polytope differs from $P$ in general, since $-\theta$ is not
usually in the $W$-orbit of $\theta$.
\end{remark}
\section{Stable coefficients}

The coefficients $a_{\lambda}$ of $N (\x; \ell )$ attached to the
vertices of $\Pi$ are called the \emph{stable coefficients}.  The goal
of this section is to show that the functional equations
\eqref{eq:coeffs}, together with the initial condition $a_{0} = 1$,
imply that the stable coefficients of $N (\x ; \ell)$ are given by the
formul\ae\ from \cite{wmd2, wmd3, wmd4}.

We caution the reader on three points.  The first is that in the
comparison that follows, it is convenient to use a slightly different
labelling convention for the vertices of $\Pi$: the element $w\in W$
now corresponds to the vertex $\theta -w^{-1}\theta $.

The second is that we are using slightly different normalizations for
Gauss sums than those found in \cite{wmd2, wmd3, wmd4}.  Hence the
formula in Theorem \ref{thm:wmd2} differs slightly from the formul\ae\
in \cite{wmd2, wmd4}, in that each factor includes a $q$-power
denominator.

Finally, in \cite{wmd4}, the twisting parameter $\ell = (l_{1},\dotsc
,l_{r})$ corresponds to the dominant weight $\theta ' = \sum
l_{i}\varpi_{i}$, whereas for us $\ell$ is attached to the regular
dominant weight $\theta = \rho +\sum l_{i}\varpi_{i}$.  Hence the
results of \cite{wmd2, wmd4} are expressed in terms of the generalized
height function $d'_{\theta}\colon \Lambda \rightarrow \Z$ defined by
$d'_{\theta }(\lambda) = \innprod{\theta+\rho }{\lambda}$.  Theorem
\ref{thm:comparison}, on the other hand, uses the function $d_\theta
(\lambda) = \innprod{\theta}{\lambda}$.  Note that $d'_{\theta '}
\equiv d_{\theta}$, so our statement of Theorem \ref{thm:wmd2} is
consistent with \cite{wmd4}.

\begin{theorem}\label{thm:wmd2}
\cite{wmd2, wmd4} Let $\lambda =\theta  -w^{-1}\theta $.  Let $A_{\lambda}$ be the
stable coefficient attached to $\lambda$ in \cite{wmd4}.  Then if
$\Phi$ is simply-laced, we have  
\begin{equation}\label{eq:gauss}
A_{\lambda} = \prod_{\alpha \in \Phi (w^{-1})} g_{1} (p^{d_{\theta }
(\alpha)-1}, p^{d_{\theta } (\alpha)})/q^{d_{\theta } (\alpha)/2},
\end{equation}
where $g_{1} (p^{a}, p^{b})$ is the quadratic Gauss sum, and where
$d_{\theta}$ is the function on the root lattice defined by
$d_{\theta} (\lambda) = \innprod{\theta}{\lambda}$.
\end{theorem}

We refer the reader to \cite{wmd2} for a precise definition of the
Gauss sum.  For our purposes all we will need to know is
\eqref{eqn:gse} below.

We now show that the stable terms of $N (\x; \ell )$ coincide with those
given by Theorem \ref{thm:wmd2}.

\begin{theorem}\label{thm:comparison}
Suppose $N (\x; \ell ) = \sum_{\lambda} a_{\lambda}\x^{\lambda}$ where
$a_{0}=1$.  Then if $\lambda =\theta -w^{-1}\theta $, the coefficient
$a_{\lambda}$ is given by \eqref{eq:gauss}.  In other words,
$a_{\lambda} = A_{\lambda}$.
\end{theorem}

\begin{proof}
We prove the theorem by induction on the length of $w$.  To begin,
note that $A_{0} = a_{0} = 1$.

Now suppose $l(\sigma_{i}w^{-1}) = l (w^{-1}) + 1$, and that the
coefficients agree on all weights attached to $u\in W$ with $l (u)\leq l (w)$.  Let $\mu = \theta
-\sigma_{i}w^{-1}\theta $.  Then $\mu = \sigma_{i}\lambda +
(l_{i}+1)\alpha_{i}$ and $\mu > \lambda$; hence we can apply
\eqref{eq:coeffs.a}, \eqref{eq:coeffs.b}, which yields
\begin{equation}\label{eq:oureqn}
a_{\mu } = a_{\lambda}\times \begin{cases}
q^{d(\sigma_{i}\lambda -\lambda +l_{i}\alpha_{i})/2}&\text{if $\varphi_{i}
(\mu )$ is even, and}\\
-q^{d (\sigma_{i}\lambda -\lambda+l_{i}\alpha_{i})/2-1/2}&\text{if $\varphi_{i}
(\mu)$ is odd}.
\end{cases}
\end{equation}
Note that the other coefficients of $N (\x; \ell )$ appearing in
\eqref{eq:coeffs.a}, \eqref{eq:coeffs.b} vanish by Theorem \ref{thm:convex} and Lemma
\ref{lem:string}.  

Applying \eqref{eqn:negroots2} in \eqref{eq:gauss}, we find
\begin{subequations}
\begin{align}
A_{\mu } &= \prod_{\alpha \in \Phi (\sigma_{i}w^{-1})} g_{1} (p^{d_{\theta }
(\alpha)-1}, p^{d_{\theta } (\alpha)})/q^{d_{\theta } (\alpha)/2}\label{eq:line1}\\
&=g_{1}(p^{d_{\theta } (w\alpha_{i})-1}, p^{d_{\theta }
(w\alpha_{i})})/q^{{d_{\theta } (w\alpha_{i})/2}} \times \prod_{\alpha
\in \Phi (w^{-1})} (g_{1} (p^{d_{\theta } (\alpha)-1}, p^{d_{\theta } (\alpha)})/q^{d_{\theta }
(\alpha)/2}).\label{eq:line2}
\end{align}
\end{subequations}
The product on \eqref{eq:line2} is just $A_{\lambda}$, which equals
$a_{\lambda}$ by induction.  Hence we must investigate the first term,
which is the Gauss sum attached to $w\alpha_{i}$.  For the quadratic
Gauss sum $g_{1} (p^{b-1}, p^{b})$ we have 
\begin{equation}\label{eqn:gse}
g_{1} (p^{b-1}, p^{b}) = \begin{cases}
q^{b-1/2}&\text{if $b$ is odd, and}\\
-q^{b-1}&\text{if $b$ is even.}
\end{cases}
\end{equation}
Hence
\begin{equation}\label{eq:theireqn}
A_{\mu} = A_{\lambda} \times \begin{cases}
q^{d_{\theta } (w\alpha_{i})/2-1/2}&\text{if $d_{\theta }
(w\alpha_{i})$ is odd, and}\\
-q^{d_{\theta } (w\alpha_{i})/2-1}&\text{if $d_{\theta }
(w\alpha_{i})$ is even.}
\end{cases}
\end{equation}
It follows that to show that \eqref{eq:oureqn} and
\eqref{eq:theireqn} agree, we must show that $d_{\theta } (w\alpha_{i})$ is even if and only if $\varphi_{i}
(\mu)$ is odd, and that $d_{\theta } (w\alpha_{i}) = d (\sigma_{i}\lambda
-\lambda +l_{i}\alpha_{i}) + 1$.

We begin by computing $d (\sigma_{i}\lambda
-\lambda +l_{i}\alpha_{i})$.  It is easy to see that 
\[
\sigma_{i}\lambda
-\lambda +l_{i}\alpha_{i} = \mu -\lambda -\alpha_{i},
\]
which implies
\[
d (\sigma_{i}\lambda
-\lambda +l_{i}\alpha_{i}) = d (\mu -\lambda)-1.
\]
Now write $\mu  = \sum k_{j}\alpha_{j}$. Then 
\[
\lambda =\sigma_{i}\mu+ (l_{i}+1)\alpha_{i} = \sum_{\substack{j\not \sim i\\
j\not =i}} k_{j}\alpha_{j} + \sum_{\adjacent{j}{i}} k_{j}
(\alpha_{i}+\alpha_{j}) + (l_{i}+1 - k_{i})\alpha_{i}.
\]
Thus 
\[
\mu -\lambda = \bigl( -\sum_{\adjacent{j}{i}}k_{j}+2k_{i}-l_{i} - 1\bigr)\alpha_{i},
\]
and 
\begin{equation}\label{eq:firstd}
d (\mu -\lambda)-1 = 2k_{i}-\varphi_{i} (\mu ) -2
\end{equation}
is our expression for $d (\sigma_{i}\lambda
-\lambda +l_{i}\alpha_{i})$.

On the other hand,
\begin{align*}
d_{\theta }(w\alpha_{i} )&= \innprod{\theta }{w\alpha_{i}}\\
	       &= \innprod{\sigma_{i}w^{-1}\theta}{-\alpha_{i}}\\
  	       &= \innprod{\theta  -\sum k_{j}\alpha_{j}}{-\alpha_{i}}\\
               &= -l_{i}-1 + \innprod{\sum k_{j}\alpha_{j}}{\alpha_{i}}.
\end{align*}
The last equation of the above gives 
\begin{equation}\label{eq:secondd}
d_{\theta }(w\alpha_{i} ) = 2k_{i}-\varphi_{i} (\mu) -1.
\end{equation}

From \eqref{eq:firstd} and \eqref{eq:secondd}, we see that $d_{\theta
} (w\alpha_{i}) = d (\sigma_{i}\lambda -\lambda + l_{i}\alpha_{i}) +
1$.  Moreover \eqref{eq:secondd} shows that the parity of $d_{\theta }
(w\alpha_{i})$ is the opposite of that of $\varphi_{i} (\beta)$.  This
completes the proof.
\end{proof}

\section{Unstable coefficients}

In this section we investigate the implications of the relations
\eqref{eq:coeffs.a}, \eqref{eq:coeffs.b} on the coefficients
$a_{\lambda}$ attached to weights other than the vertices of $\Pi$.
Such coefficients are called \emph{unstable coefficients} \cite{wmd2,
wmd3, wmd4}.  The main result of this section, Theorem
\ref{thm:uniquelydetermined}, is that in the untwisted case $\theta =
\rho$, the numerator polynomial and hence $f (\x)$ is uniquely
determined by \eqref{eq:coeffs.a}, \eqref{eq:coeffs.b} and the
normalization condition $a_{0} = 1$.  We conclude by discussing the
extent to which $N (\x ;\ell)$ is not uniquely determined.

Recall that a weight $\lambda $ is \emph{regular} if it lies in the
interior of a Weyl chamber.  Equivalently, the stabilizer $P= P (\lambda )$
of $\lambda $ in $W$ is trivial.  The stabilizer $P$ of any 
weight is a subgroup of $W$ generated by a subset of the simple
reflections.  Indeed, if $\lambda $ lies on the hyperplane fixed by
$\sigma_{i}$, then $\sigma_{i}\in P$, and such simple reflections
generate $P$.  Any subgroup generated by a subset of the simple
reflections is called a \emph{standard parabolic subgroup}.  We recall
the following basic fact about such subgroups, whose proof can be
found, for example, in \cite{humph}:

\begin{prp}\label{prop:cosets}
Let $P\subset W$ be a standard parabolic subgroup.  Then any coset
$wP$ contains a unique element $w^{P}$ of maximal length.
\end{prp}

Now let $N (\x; \ell ) = \sum a_{\lambda} \x^{\lambda}$ be our polynomial.
By Theorem \ref{thm:convex} the support of $N (\x; \ell )$ consists at most
of the monomials $\x ^\lambda$ where $\lambda \in \Lambda$ lies in the
convex hull of the point $\theta  -w\theta , w\in W$.  Such $\lambda$
correspond to the weights of the representation $V_{\theta }$ with
highest weight $\theta $, after shifting.  The precise connection is as
follows.  Let $\Theta$ be the set of \emph{dominant weights} of
$V_{\theta }$.  Then the support of $N (\x;\ell )$ consists of all monomials
$\x^{\lambda }$ where $\lambda$ has the form
\[
\lambda = \theta - w\xi , \quad w\in W,\, \xi \in \Theta.
\]

For $\xi \in \Theta$, let $O_{\xi }$ be the $W$-orbit $\{ \theta -
w\xi\mid w\in W \}$, and let $\sO = \{O_{\xi} \mid \xi \in \Theta \}$.
The set $\Theta$ is naturally a poset by the usual partial order on
weights (note that $\theta $ is the maximal element), and we use this
to turn $\sO$ into a poset: $O_{\xi }\leq O_{\xi '}$ if and only if
$\xi \leq \xi '$.

For any $\xi \in \Theta$, let $\Pi_{\xi }$ be the convex hull in
$\Lambda_{\R}$ of the points in $O_{\xi }$.  If $\xi $ is regular
then $\Pi_{\xi }$ is isomorphic to a permutahedron, with vertices in bijection
with $W$.  Otherwise, the orbit $O_{\xi }$ has fewer than $|W|$
vertices.  Indeed, if $P = P (\theta)\subset W$ is the stabilizer of
$\xi $, then $|O_{\xi }|$ equals the number of cosets $|W/P|$.
Accordingly, 
we can also write
\[
O_{\theta} = \{ \theta -  wP\xi  \mid w\in W\}.
\]

\begin{remark}
Although we do not need it, it is known that $\Pi_{\xi }$ is a
degeneration of $\Pi_{\theta }$ obtained by contracting the edges in
$\Pi_{\theta }$ labelled by the simple reflections in $P$.  Another
way to express the relationship between $\Pi_{\xi}$ and $\Pi_{\theta}$
is to observe that $\Pi_{\xi }$ is obtained from $\Pi_{\theta }$ by
parallel translation of the facets of $\Pi_{\theta }$ until some of
the faces collapse.  This has the consequence that essentially the
same system of inequalities \eqref{eq:inequalities} describes the
polytope $\Pi_{\xi }$.  In particular, $\Pi_{\xi }$ is described by
\begin{equation}\label{eq:inequalities2}
\innprod{w\weight_{k}}{\x - (\theta  -wP\xi) }\geq 0, \quad k=1,\dotsc ,r,
\quad w\in W.
\end{equation}
As before this system is redundant.  Also, some inequalities do
not define facets, and instead are active only on higher-codimension faces.
Nevertheless $\Pi_{\xi}$ is cut out by the system \eqref{eq:inequalities2}.
\end{remark}

We are now ready to prove the geometric results that allow us to
analyze unstable coefficients.  The main point is the following
generalization of Lemma \ref{lem:string}:

\begin{lemma}\label{lem:string.distinct}
Let $\lambda = \theta  -u P\xi $ be a vertex of $\Pi_{\xi }$, where
$P=P (\xi)$, and suppose $u\in
uP$ is the unique maximal element in this coset.  Let $w=\sigma_{k}u$
and suppose $l (w)>l (u)$.  Then $\mu =\theta  - wP\xi $ is a different
vertex of $\Pi_{\xi }$.  Moreover, if any point of the form
$\mu +m\alpha_{k}$, $m\geq 1$ lies in an orbit $O\in
\sO$, we have $O > O_{\xi }$.  Similarly, if any point of
the form $\lambda - m\alpha_{k}$, $m\geq 1$ lies in an orbit
$O\in \sO$, we have $O > O_{\xi }$.
\end{lemma}

Before we prove Lemma \ref{lem:string.distinct}, we need another lemma
about the geometry of dominant weights:

\begin{lemma}\label{lem:dominance}
Let $\xi$ and $\eta$ be weights such that $\eta>\xi $.  Let $\eta '$
be the unique dominant weight in the $W$-orbit of $\eta$.  Then $\eta
'>\xi$.
\end{lemma}

\begin{proof}
Let $C (\Phi^{+})$ be the cone generated by the positive roots.  Then
$\eta >\xi$ implies $\eta -\xi \in C (\Phi^{+})$.  On the other hand
if $\eta '$ is the dominant weight in the orbit of $\eta$, then
certainly $\eta '-\eta \in C (\Phi^{+})$.  But then $\eta '-\xi \in C
(\Phi^{+})$, since $C (\Phi^{+})$ is convex.  This completes the
proof.
\end{proof}

\begin{proofof}{ of Lemma \ref{lem:string.distinct}}
First, it is clear that $\lambda \not =\mu $, since
$l (w)>l (u)$ and $u$ is the maximal element of the coset $uP$.

We now show $O>O_{\xi}$, where $O$ is the orbit corresponding to $\mu
+m\alpha_{k}$, $m>1$.  We have $\xi =w^{-1}(\theta - \mu )$.  Let
\[
\eta  = w^{-1}(\theta - \mu -m\alpha_{k}) = \xi -mw^{-1}\alpha_{k}, 
\]
and let $\eta '$ be the unique dominant weight in the $W$-orbit of
$\eta$.  Then $O=O_{\eta  '}$.  Now 
\begin{equation}\label{eq:negative}
\Phi (w^{-1}) = \Phi (u^{-1}\sigma_{k}) = \sigma_{k}\Phi (u^{-1}) \cup
\{\alpha_{k} \},
\end{equation}
which implies
\[
w^{-1}\alpha_{k}\in \Phi^{-}.
\]
Thus $\eta > \xi$.  By Lemma \ref{lem:dominance} we have $\eta '>\xi$,
which proves $O>O_{\xi}$.

The statement about $\lambda - m\alpha_{k}$ is proved in almost the
same way.  The computation boils down to 
\[
\alpha_{k}\not \in \Phi (u^{-1}).  
\]
This is clearly true, since by \eqref{eq:negative} the set
$\sigma_{k}\Phi (u^{-1})$ consists of positive roots.
\end{proofof}

Hence in applications of \eqref{eq:coeffs.a}, \eqref{eq:coeffs.b}, if
$\lambda$ and $\mu$ are weights with $\mu =\sigma_{k}\lambda +
(l_{k}+1)\alpha_{k}$ and $\mu >\lambda $, then we know that the
weights $\lambda-2\alpha_{k}$, $\lambda -\alpha_{k}$, $\mu
+\alpha_{k}, \mu +2\alpha_{k}$ live in bigger orbits in $\sO$ and are
attached to previously determined coefficients.

We likewise have a version of Lemma \ref{lem:string} even if
$\sigma_{k}$ fixes a vertex:

\begin{lemma}\label{lem:string.same}
Let $\lambda = \theta -uP\xi $ be a vertex of $\Pi_{\xi }$. Let
$w=\sigma_{k}u$ with $l (w)>l (u)$, and suppose $\mu = \theta  - wP\xi$
equals $\lambda $. 
Then if any point of the form $\mu + m\alpha_{k}$,
$m\geq 1$ lies in an orbit $O\in \sO$, we have $O> O_{\xi
}$.  Similarly if any point of the form $\lambda - m\alpha_{k}$,
$m\geq 1$ lies in an orbit $O$, we have $O> O_{\xi
}$.  
\end{lemma}

\begin{proof}
The proofs of both statements are essentially the same as those of
Lemma \ref{lem:string.distinct}, even though the points $\lambda$,
$\mu$ coincide.  Again the key points are that $\alpha_{k}\not \in
\Phi (u^{-1}), \alpha_{k}\in \Phi (w^{-1})$.
\end{proof}

\begin{remark}
It is perhaps inaccurate to describe Lemmas \ref{lem:string.distinct} and
\ref{lem:string.same} as generalizations of Lemma \ref{lem:string},
since the statements are so different.  However they really are the
same, since it is the same geometric fact about left descent sets that
is behind all of them.

Also, the inequalities \eqref{eq:inequalities2} are lurking here as
well.  The points $\mu + m\alpha_{k}$, $\lambda -m\alpha_{k}$ lie
outside $\Pi_{\xi }$.  The fact that these points violate
\eqref{eq:inequalities2} again boils down to $\alpha_{k}\not \in \Phi
(u^{-1}), \alpha_{k}\in \Phi (w^{-1})$.
\end{remark}

We can now prove the main result of this section.

\begin{theorem}\label{thm:uniquelydetermined}
Let $N (\x; \ell )$ be the numerator, normalized so that $a_{0}=1$.
Suppose that $\theta$ is the only regular dominant weight in the
representation $V_{\theta}$ of highest weight $\theta$.  Then all
other coefficients of $N (\x; \ell )$ are uniquely determined by the
functional equations \eqref{eq:coeffs.a}, \eqref{eq:coeffs.b}.  
\end{theorem}

\begin{proof}
The proof is by descending induction over the orbit poset $\sO$.  To
begin, we know from Theorem \ref{thm:comparison} that all the
coefficients attached to the elements of the orbit $O_{\theta }$ are
uniquely determined once we know $a_{0}=1$.

Now fix an orbit $O_{\xi }$, where $\xi \in \Theta$ is different from
$\theta $, and assume we have determined the coefficients attached to
all orbits $O$ with $O>O_{\theta}$.  By assumption $\xi $ is not
regular, so we can find a simple functional equation from
\eqref{eq:coeffs.a}, \eqref{eq:coeffs.b} relating the corresponding
coefficient $a_{\lambda }$, $\lambda =\theta -\xi$, to itself.  It is trivial to see that if
$\varphi_{k} (\lambda)$ is even, then $a_{\lambda }$ appears on both
sides of \eqref{eq:coeffs.a} with different coefficients; if
$\varphi_{k} (\lambda)$ is odd then clearly $a_{\lambda }$ appears on
both sides of \eqref{eq:coeffs.b} with different coefficients. By
Lemma \ref{lem:string.same}, all other $a_{\lambda '}$ in
\eqref{eq:coeffs.a}, \eqref{eq:coeffs.b} come from previously
determined orbits.  Thus $a_{\lambda }$ is determined.

Now, successively applying Lemma \ref{lem:string.same} we can
determine the remaining coefficients of the form $\theta -w\xi$, where
$w \in W$.  The basic
strategy is the same as in the proof of Theorem \ref{thm:comparison};
if there is more than one point in this orbit, all one needs to be
able to do is move from one to another by left multiplication by a
simple reflection
\[
 \mu := \theta -
wP\xi  \longrightarrow \theta - \sigma_{k}wP\xi  =:\mu' ,
\]
where $w$ is maximal in the coset $wP$ and $l (\sigma_{k}w)>l (w)$.
Then Lemma \ref{lem:string.distinct} shows that $a_{\mu }$ and the
coefficients from higher orbits determine $a_{\mu' }$.
\end{proof}

\begin{corollary}\label{cor:untwistedunique}
The only regular dominant weight for the representation $V_\rho$ is $\rho$.
Thus the numerator in the untwisted case,
$N(\x)=N (\x ; (0,\dots,0))$, is uniquely determined by the
functional equations \eqref{eq:coeffs.a}, \eqref{eq:coeffs.b} after
setting $a_{0}=1$.
\end{corollary}

\begin{proof}
Suppose $\lambda\not =\rho $ is another regular dominant weight for this
representation.  Write $\lambda = \sum c_{i}\weight_{i}$.  Then we
must have each $c_{i}\not =0$, since $\lambda$ is regular.   Since
$\rho = \sum_{i}\weight_{i}$, it follows that 
$\rho -\lambda $ must be a linear combination of the $\weight_{i}$
with nonpositive coefficients.

On the other hand, $\rho -\lambda$ is a nonnegative linear combination
of the simple roots, since $\rho$ is higher than $\lambda$ in the
partial order.  The fundamental weights are themselves positive
rational linear combinations of the simple roots, as one sees by
examining the inverse Cartan matrix for any simple complex Lie
algebra.  Thus $\rho -\lambda$ is simultaneously a nonpositive and a
nonnegative linear combination of the fundamental weights.  This means
$\rho =\lambda$, a contradiction.  Hence $\rho$ is the only regular
dominant weight in the representation $V_{\rho}$, and by Theorem
\ref{thm:uniquelydetermined} the polynomial $N (\x)$ is uniquely
determined. 
\end{proof}

\begin{remark}
There are other regular dominant weights $\theta\not =\rho$ such that
$V_{\theta}$ satisfies the conditions of Theorem
\ref{thm:uniquelydetermined}.  For example, computations show that for
$\Phi = A_{r}$ with $r\leq 5$ the representation $V_{\rho +
\varpi_{1}}$ has a unique regular dominant weight; presumably this
representation does for all $r$.  We do not know another
characterization of weights for any given $\Phi$ with this property.
\end{remark}

\begin{remark}
One can also ask how big the space of possible numerators $N (\x ;
\ell)$ can be if $\theta$ does not satisfy the conditions of Theorem
\ref{thm:uniquelydetermined}.  Examples for $A_{2}$ show that the
complex dimension of this space apparently equals the number of
regular dominant weights of the representation $V_{\theta}$.  This
should be true for all $\Phi$, although we have not checked the
details.
\end{remark}


\section{The global multiple Dirichlet series}\label{s:twmds}
In this final section, we describe precisely how a multiple Dirichlet
series can be built up out of the $p$-parts.  We follow the methods of
\cite{qmds}, where the untwisted case was worked out in detail.  We
remark, however that the $H$ function  defined below is the analogue
of the $H$ function of \cite{wmd1, wmd2,wmd3, wmd4} rather than that
of \cite{qmds}.  For the relation between the two, see Remark 4.3 of
\cite {qmds}.  
    
For simplicity, we work over $\mathbb{Q}, $ the field of rational
numbers, and let $\prs{d}{m}$ denote the usual quadratic residue
symbol for $d,m$ odd and relatively prime.  For an arbitrary global
field $K,$ one needs to work over the ring ${\mathcal{O}}_S$ of
$S$-integers of $K$ for a sufficiently large set of primes $S$ and to
replace $p$ by $q=|{\mathcal O}_S/p{\mathcal O}_S|$.  Additionally,
some care needs to be taken to define the quadratic residue symbol
properly.  We refer the reader to \cite{ff} and \cite{qmds} for
details.

Given an odd prime $p$ and twisting parameter $\ell,$  we write
\begin{equation}
N(\x;\ell)=\sum_{\lambda \in \Lambda}
a_{\lambda}(p,\ell)\x^{\lambda }.
\end{equation}

Fix an $r$-tuple of positive odd integers $\ttt=(t_1, t_2,\ldots,
t_r).$ The goal of this section is to  define
the $\ttt$-twisted multiple Dirichlet series 
associated to the root system $\Phi$ for $n=2.$  Actually, we will
need to introduce a  family of series 
$$Z(s_1, \ldots, s_r;\Psi, \ttt;\Phi)$$
where $\Psi$ ranges over $r$-tuples
$\Psi=(\psi_1,\psi_2,\ldots,\psi_r)$
of Dirichlet characters unramified outside of $2$.  We abbreviate this
series by $Z(\s;\ttt, \Psi),$ and it is understood that  $\Phi$
remains fixed.  Each series will be a sum over $r$-tuples of odd positive
integers.

We call an $r$-tuple ${\bf m}=(m_1,\ldots,m_r)$ of positive
integers {\em odd} if each of the $m_i$'s is odd.  We denote
by $\Psi({\bf m})$ the product 
$$\prod_{i}\psi_i(m_i)$$  
and by $H({\bf m};\ttt)$ the coefficient $H(m_1,m_2,\ldots,m_r;\ttt)$ 
defined below.
\begin{defin} \label{def:H} The coefficient $H(m_1,m_2\ldots,m_r;\ttt)$ is
  defined by the following two conditions:
\begin{enumerate}
\item Suppose  $\mathbf{m} = (p^{k_{1}},\dotsc , p^{k_{r}})$, where $p$
is an odd prime and $p^{l_i}||t_i.$
Suppose $\lambda =\sum_{i=1}^r k_i\alpha_i\in
\Lambda.$  Then
\begin{equation*}H(p^{k_1},\ldots,p^{k_r};\ttt)=
a_\lambda(p;\ell) \end{equation*}
where $\ell=(l_1,\ldots, l_r).$ 
\item  Given  $m_j, m_j^\prime$ odd with $(m_1m_2\cdots
  m_r, m_1^\prime m_2^\prime\cdots m_r^\prime)=1$ we have 
\begin{equation*}
\frac{H(m_1m_1^\prime, \ldots, m_rm_r^\prime;\ttt )}
{H(m_1, \ldots, m_r;\ttt)H(m_1^\prime, \ldots, m_r^\prime;\ttt )}=
\prod_{\substack{i,j \adj\\i<j}} \prs{m_i}{m_j^\prime}
\prs{m_i^\prime} {m_j}
\end{equation*}
\end{enumerate}
\end{defin}

Finally, for an $r$-tuple $\s=(s_1,\ldots,s_r)$ of
complex numbers, define  
\begin{equation}\label{def:Ztwisted}
Z(\s;\ttt,\Psi)=N(\s)\sum_{{\bf m}=(m_1,m_2,\ldots, m_r)\ odd}
\frac{\Psi({\bf m})H({\bf m};\ttt)}
{\prod_j m_j^{s_j}} \prod_{i=1}^{r} \prs{t_i^\#}{\hat m_i}
\end{equation}
where $t_i^\#$ is the squarefree part of $t_i$,  $\hat m_i$ is the 
part of $m_i$ relatively prime to $t_i^\#$ 
and $N(\s)$ is the
normalizing zeta factor
\begin{equation}\label{eqn:nzf}
N(\s)=\prod_{\alpha\in\Phi^+} \zeta(2\langle \alpha, {\bf
s}\rangle-d(\alpha)+1), \ \ \langle\alpha,{\bf s}\rangle=\alpha_1
s_1+\cdots+\alpha_r s_r. 
\end{equation}
The series \eqref{def:Ztwisted} converges for $\Re(s_i)\gg 1$, $1\leq
i\leq r$.

Now, given the invariance of the $p$-parts of $Z(\s;\ttt, \Psi)$ under
the action of $W$ defined by \eqref{wiactionTwisted}, we may mimic the
techniques of Section 5 of \cite{qmds} to show the following: for
fixed twisting parameter $\ttt$, the vector of all $Z(\s;\ttt,\Psi)$
has analytic continuation and satisfies functional equations relating
the values at $\s=(s_1,\ldots, s_r)$ to the values at
$\sigma_{j_0}\s=(s_1^\prime,\ldots,s_r^\prime)$ for $j_0=1,2,\ldots,
r,$ where
\begin{equation}\label{eqn:actionOnS}
s_j^\prime=\left\{\begin{array}{ll}
s_j+s_{j_0}-1/2 & \mbox{\ if $j$ and $j_0$ are adjacent,}\\
1-s_{j_0} & \mbox{\ if $j=j_0$, and}\\
s_j & \mbox{\ otherwise.}
\end{array}\right.
\end{equation}
These functional equations are involutions generating a group of
functional equations of $Z(\s;\ttt, \Psi).$ We then deduce

\begin{theorem}\label{thm:main}
Fix a twisting parameter $\ttt$.  Each function $Z(\s;\ttt,\Psi)$ has
analytic continuation to $\C^r$.  The collection of these functions as
$\Psi$ ranges over $r$-tuples of Dirichlet characters unramified
outside of $2$ satisfies a group of functional equations 
isomorphic to $W$.  Finally, each $Z(\s;\ttt,\Psi)$ is analytic
outside the hyperplanes $(w\s)_j=1$ for $w\in W, 1\leq j\leq r$, where
$(w\s)_j$ denotes the $j^{th}$ component of $w\s$.
\end{theorem}

The proof of the theorem is very similar to the proof of Theorem 5.5 of
\cite{qmds}, and so we leave the details to the reader.

\bibliographystyle{amsplain_initials}
\bibliography{ratfun}

\end{document}